\documentclass{llncs}

\usepackage[utf8]{inputenc}
\usepackage{mystyle,amssymb,dsfont,amsmath,mathtools}

\DeclareMathOperator{\KL}{KL}

\DeclareMathOperator{\LSE}{LSE} 

\def\ones{\mathds{1}}

\newcommand{\Loss}{\mathcal{L}}

\newcommand{\mysubsection}[1]{\vspace{1mm}\noindent\textbf{#1.}~}
\newcommand{\myparagraph}[1]{\noindent\textit{#1.}~}

\begin{document}
\frontmatter          
%
%
%
\title{
Optimal Transport for Diffeomorphic~Registration
}
\titlerunning{Optimal transport for diffeomorphic registration}  

%




\author{Jean Feydy\inst{1,2} \and Benjamin Charlier\inst{3,5} \and \\  François-Xavier Vialard\inst{4,6} \and Gabriel Peyré\inst{1,5} }
 
\institute{DMA -- École Normale Supérieure, Paris, France
\\ \email{jean.feydy@ens.fr}, \email{gabriel.peyre@ens.fr}
\and
CMLA -- ENS Cachan, Cachan, France
\and
Institut Montpelliérain Alexander Grothendieck, \\ Univ. Montpellier, Montpellier, France \\
\email{benjamin.charlier@umontpellier.fr}
\and
Univ. Paris-Dauphine - PSL Research, Paris, France \\ \email{vialard@ceremade.dauphine.fr} 
\and
CNRS, Paris, France
\and INRIA Mokaplan, Paris, France
}

 \authorrunning{-}   
 \tocauthor{-}

\maketitle              

\begin{abstract}
This paper introduces the use of unbalanced optimal transport methods as a similarity measure for diffeomorphic matching of imaging data. 
The similarity measure is a key object in diffeomorphic registration methods that, together with the regularization on the deformation, defines the optimal deformation. Most often, these similarity measures are local or non local but simple enough to be computationally fast.
We build on recent theoretical and numerical advances in optimal transport to propose fast and global similarity measures that can be used on surfaces or volumetric imaging data.
This new similarity measure is computed using a fast generalized Sinkhorn algorithm. We apply this new metric in the LDDMM framework on synthetic and real data, fibres bundles and surfaces and show that better matching results are obtained.
\end{abstract}


\section{Introduction}


State of the art methods in deformable registration of medical image data are often variational methods that estimate an optimal deformation by minimizing an energy which is the sum of two terms: a data fidelity term which quantifies how well the source data has been aligned with the target data; a regularization term on the deformation, which is necessary to make the registration problem well-posed, see \cite{Sotiras2013} for a detailed overview.
In the following, we will be interested in diffeomorphic image matching, resulting in a non-convex optimization problem and in estimation of local minima of the energy. The choice of the data fidelity is thus crucial to avoid poor local minima and to obtain meaningful registration. In this paper we propose a ``global'' fidelity making use of optimal transport theory and fast entropic approximation schemes. 

\mysubsection{Previous Works}
\myparagraph{Data fidelity for registration}
Several similarity measures have been introduced in the literature which emphasizes its crucial role.
For instance, for dense image registration, when only the intensity information is considered, the sum of squared differences (SSD)~\cite{Gee1993} or generalized versions of the cross-correlation such as normalized cross-correlation \cite{Avants2008} are used which correct for possible intensity biases.
When the images have already been segmented, a lot of interest has been devoted to matching between shape using parametrization invariant data fidelities. Let us mention currents, varifolds and more recently functional currents. 
All these shape metrics can be understood as a non-local SSD. Most often, the similarity measures are local or at most non-local in order to keep the computational cost low. 

\myparagraph{Unbalanced optimal transport}
The basic theory of optimal transport (OT) (see~\cite{santambrogio2015optimal}) defines a distance between probability distributions using the amount of effort needed to transport mass from one measure to another one. These distances are appealing for geometric problems such as shape registration because they are sensitive to spatial displacements of the shape. A recent breakthrough has been made in extending OT to the case of measures of different total mass \cite{LieroMielkeSavareLong,2016-chizat-focm}, which is crucial for practical applications where mass can vary because of changes in shape scales, or because of mass creation/destruction processes.
%
%
%
%
The use of OT as a fidelity term is not new, see~\cite{montavon2016wasserstein} for an example in machine learning. Quite surprisingly, to the best of our knowledge, it has never been used for diffeomorphic registration purpose. 

\myparagraph{Entropic regularization}
A critical aspect of OT is that it was involved numerically. 
This situation has however radically changed in the last few years, thanks to the introduction of an efficient entropic approximation scheme~\cite{cuturi-2013}, 
which is (i) efficient and easily parallelizable, (ii) and leads to a smooth differentiable data fidelity term~\cite{CuturiDoucet}. 
This scheme also applies to unbalanced OT~\cite{2016-chizat-sinkhorn}.

\mysubsection{Contributions} 
This paper proposes a new non-local geometric similarity measure based on the recently developed theory of unbalanced OT and fast entropic solvers.
The resulting hybrid pipeline is able to combine the strength of both diffeomorphic models presented in Sec.~\ref{SecDiffReg} and entropic unbalanced OT (convex optimization with order of magnitude faster solvers) detailed in Sec.~\ref{Sec:Transport}. As shown by simulations on synthetic and real data in Sec.~\ref{sec-numerics}, OT can thus seamlessly be integrated in state-of-the-art registration methods, enabling a long-range non-local attraction term toward the target (which is crucial to avoid poor local minima) while remaining able to match intricate fine scale details. 
The numerical code to reproduce the figures of this article is provided as supplementary material. 
\section{Diffeomorphic Registration} \label{SecDiffReg}

\mysubsection{Data Representation with Measures}
\label{sec:measures_repr}
After a pre-processing step, it is often possible to efficiently represent medical image data as measures over a conveniently chosen space  $\XX$. This representation of data is essentially motivated by methodological and numerical purposes to help the design of fidelity terms. 
%

For registration of segmented shapes (curves, surfaces), we advocate the use of a lifted
features space $\XX = \RR^d \times \SS^{d-1}$: each segment (\emph{respectively} triangle)
of the curve (\emph{resp.} surface)
can be represented as a Dirac of mass $p_i$ equal to its length (\emph{resp.} area), placed at location $(a_i, u_i) \in \XX$, where $a_i$ is the center and $u_i$ the unit-length normal of the shape element.

Eventually, source and target are thus represented as 
\eql{\label{eq-discr-measured}
	\textstyle \mu = \sum_{i \in I} p_i \de_{x_i}
	\qandq
	\nu = \sum_{j \in J} q_j \de_{y_j}
}
where $x_i, y_j \in \XX$ are the sampled features, $p_i, q_j \geq 0$ are the associated masses, and $\de_x$ denotes the Dirac mass at some point $x \in \XX$. 

\mysubsection{Diffeomorphic Registration and Data Fidelity}
Variational diffeomorphic registration of shape $\mu$ onto the shape $\nu$ consists in the minimization of 
\begin{equation}
\mathcal{E}(f) \eqdef \mathcal{R}(f) + \Loss(f_\ast\mu,\nu)
\end{equation}
where $\mathcal{R}$ is the regularization term on the diffeomorphism $f : \RR^d \rightarrow \RR^d$ and $\Loss$ is the data fidelity.
The notation $f_* \mu$ stands for the data $\mu$ deformed by $f$.  

The most simple data fidelity terms are derived from Euclidean norms using a smoothing operation against a kernel $G_\si$ of width $\si$
\eql{\label{eq-rkhs-loss}
	\Loss(\mu,\nu)=\int_\XX \Big( 
	    \int_\XX G_\si(x,y) (\d\mu(x)-\d\nu(x))
	    \Big)^2 \d y.
} 
This class of losses has been used extensively for shape matching (see for instance~\cite{Vaillant2005,CharonT13}) and is also popular in machine learning under the name Kernel mean embedding~\cite{sejdinovic2013}.
Its limitations for diffeomorphic registration are discussed in section~\ref{sec:protozoa}.

\mysubsection{Optimization Scheme}
In all the cases of interest, the action $f_* \mu$ of $f$ on a finite discrete measure of the form~\eqref{eq-discr-measured} can be written using a finite dimensional vector $\th$ of parameters. 
This formulation includes non-parametric methods since the parametrization may depend on the input shape, as in LDDMM methods.
We write down this action as
$
	f_* \mu = \mu_\th \eqdef \sum_i p(\th)_i \de_{x(\th)_i}. 
$

Registration is achieved by computing a local minimizer of
\eql{\label{eq-min-register}
	\umin{\th} \mathcal{E}(\th) = \mathcal{R}(\th) +  \Loss(\mu_\th,\nu)
}
using a descent-based method, typically initialized for $\th=0$ using $\mu_0=\mu$, 
i.e. $p(0)_i=p_i, x(0)_i=x_i$. 
%
%
From a computational perspective, we simply need to provide the gradient of the functional which reads, thanks to the chain rule
\eql{\label{eq-chain-rule}
	\nabla \mathcal{E}(\th) = \nabla \mathcal{R}(\theta) + [\partial x(\th)]^*( \nabla_x \Loss(\mu_\th,\nu) )
		+ [\partial p(\th)]^*( \nabla_p \Loss(\mu_\th,\nu) ),
}
where $([\partial x(\th)]^*, [\partial p(\th)]^*)$ are the adjoints of the Jacobians of the maps $(\th \mapsto x(\th),\th \mapsto p(\th))$ and $(\nabla_x \Loss(\mu,\nu),\nabla_p \Loss(\mu,\nu))$ are the gradients of the maps $x \mapsto \Loss(\sum_i p_i \de_{x_i},\nu)$ and $p \mapsto \Loss(\sum_i p_i \de_{x_i},\nu)$. 

\section{Optimal Transport Data Fidelity}\label{Sec:Transport}

This section proposes a new class of data fidelity $\Loss(\mu,\nu)$ between measures, using the recently proposed framework of unbalanced optimal transport between positive measures. 
%
%

\mysubsection{Unbalanced Regularized Optimal Transport}
We consider two input discrete measures~~\eqref{eq-discr-measured}. 
In OT, the transportation can be described by a joint distribution defined on $\XX \times \XX$ coupling the two measures.
For discrete inputs, this coupling is an array 
$\ga=(\ga_{i,j})_{i,j}=\text{``displaced mass from $x_i$ to $y_j$''}$
of positive numbers, whose marginals 
$
    \textstyle \ga \ones_J \eqdef ( \sum_j \ga_{i,j} )_i
$
and
$
    \ga^\top \ones_I \eqdef ( \sum_i \ga_{i,j} )_j
$
should be equal (for classical balanced) or close (for unbalanced)  to the input measures
$\mu$ (source) and $\nu$ (target). 

An approximate OT discrepancy is obtained by looking for an optimal coupling as
the \emph{regularized} optimal transport cost is given by
\eql{\label{eq-primal-relax}
	W_{\epsilon,\rho}(\mu,\nu) \!\eqdef\!\! \umin{\ga \in \RR_+^{I \times J}}
		\sum_{i,j} c(x_i,y_j) \ga_{i,j}  - \epsilon H(\ga)
		+ \rho \KL(\ga \ones_J|p) 
		+ \rho \KL(\ga^\top \ones_I|q).
}
Here the Kullback-Leibler divergence and entropy read
\begin{align*} 
	\KL(h|p) &\eqdef \textstyle\sum_i h_i \log\pa{ \frac{h_i}{p_i} } - h_i + p_i, \quad
	 H(\ga) &\eqdef-{\textstyle\sum_{i,j}} \ga_{i,j} (\log(\ga_{i,j})-1), 
\end{align*}
and $c(x_i,y_j)$ is the cost of displacing a \emph{unit} amount of mass between positions $x_i$ and $y_j$.
In the problem~\eqref{eq-primal-relax}, $\epsilon \geq 0$ controls the degree of regularization, and setting $\epsilon=0$, $\rho=+\infty$ recovers the usual OT. 
The influence of both parameters $\epsilon$ and $\rho$ is discussed in Section~\ref{sec-numerics}.
\mysubsection{Generalized Sinkhorn Algorithm}
Following~\cite{2016-chizat-sinkhorn}, one can check that the optimal $\ga$ can be written in the form
\eql{\label{eq-primal-dual}
	\ga = \exp(K(u,v))
	\text{ where }
	\foralls (i,j) \in I \times J, \quad
	K(u,v)_{i,j} \eqdef \frac{u_i+v_j-c(x_i,y_j)}{\epsilon} , 
}
which thus only requires the computation of two ``dual'' vectors $(u,v) \in \RR^I \times \RR^J$. In addition, these two vectors can be computed using an adaptation of the classical Sinkhorn algorithm. 
Introducing the log-sum-exp operator
$
	\LSE_I( K ) = \log \sum_{i\in I} \exp(K_{i,j}) \in \RR^J
$
(and similarly for $\LSE_J$, doing summation over $j \in J$), and starting from $(u,v)=(0_I,0_J)$, Sinkhorn's iterations read
\begin{equation}\label{eq-sinkhorn}
\begin{aligned}
	u &\leftarrow \la u
	+ \epsilon\la \log(p)
	-\epsilon\la  \LSE_J(K(u,v)) \\
	v &\leftarrow \la v 	
	+ \epsilon\la \log(q)
	-\epsilon\la  \LSE_I(K(u,v))
	\end{aligned}
\end{equation}
where we defined $\la \eqdef \frac{\rho}{\rho+\epsilon}$. The output of the algorithm is then $\ga=\exp(K(u,v))$. Note that when $\rho=+\infty$ (balanced case), $\la=1$ and these iterations correspond to a stable implementation over the log-domain of the well-known Sinkhorn algorithm (which is usually written using multiplicative updates, which is unstable for small $\epsilon$). 
This algorithm is known to converge linearly to the optimal coupling.

\mysubsection{Derivatives of the OT Fidelity} \label{part:grad}
%
%
A proof similar to the balanced case (see~\cite{CuturiDoucet}) shows that for $\epsilon>0$, $(p,x) \mapsto W_{\epsilon,\rho}(\mu,\nu)$ is smooth and its gradient reads
\begin{align}\label{eq-ot-gradients}
	\nabla_p W_{\epsilon,\rho}(\mu,\nu) = 
		\rho (1-e^{-\frac{u}{\rho}}) \qandq 
	\nabla_x W_{\epsilon,\rho}(\mu,\nu) = (\textstyle \sum_j \ga_{i,j} \partial_1 c(x_i,y_j) )_{i}
\end{align}
(with the convention $(1-e^{-\frac{u}{\rho}})=u$ for $\rho=+\infty$)
where $\ga$ is the solution of~\eqref{eq-primal-relax} and $u$ is the limit of Sinkhorn iterations~\eqref{eq-sinkhorn}. 
Here $\partial_1 c$ is the derivative of the cost $c$ with respect to the first variable. 

\section{Numerical Results}
\label{sec-numerics}

In this section, we showcase the use of unbalanced OT as a versatile fidelity term for registration.
Our code is available: \url{\texttt{github.com/jeanfeydy/lddmm-ot}}.


\mysubsection{Practical Implementation of OT Fidelity}
%
To use OT on real data, one simply needs: an appropriate mapping from the dataset to the space of measures on a features space $\XX$ ;
a cost function $c(x,y)$ on $\XX\times\XX$, with values in $\RR^+$.

In the case of curves/surfaces, following the construction of Section~\ref{sec:measures_repr},
one needs to choose a cost function on the (positions,normals) product 
$\XX=\RR^d\times\SS^{d-1}$.
One can use the canonical distance between $x=(a,u)$ and $y=(b,v)$
\begin{align}
c(x,y) &= \norm{a-b}^2 +  \alpha \,d_{\SS^{d-1}}^2(u,v), \label{eq:cprod}\\
\text{or, for instance, ~~}
c(x,y) &= \norm{a-b}^2 \cdot \Big( 1 + \alpha \, \big( 1 - \dotp{u}{v}^k \big) \Big)\label{eq:cprodbis}
\end{align}
where $\alpha\geqslant 0$ and $d_{\SS^{d-1}}$, $k$ parametrizes the angular selectivity of the registration.
Doing so, choosing $\alpha = 0$ allows one to retrieve the standard Wasserstein distance between shapes,
whereas using $d_{\SS^{d-1}}(u,v)  = (1 - \dotp{u}{v})$, $k=1$ 
(resp. $d_{\SS^{d-1}}(u,v) = (2 - 2  \dotp{u}{v}^2)$, $k$ even) can be seen as using globalized variants of the classical
currents~\cite{Vaillant2005} (resp. varifold~\cite{CharonT13}) costs.

\begin{figure}[b!]
\begin{tabular}{@{~~~~~~~~~~~~}c@{}c@{}c@{}}	
\includegraphics[width=.22\linewidth]{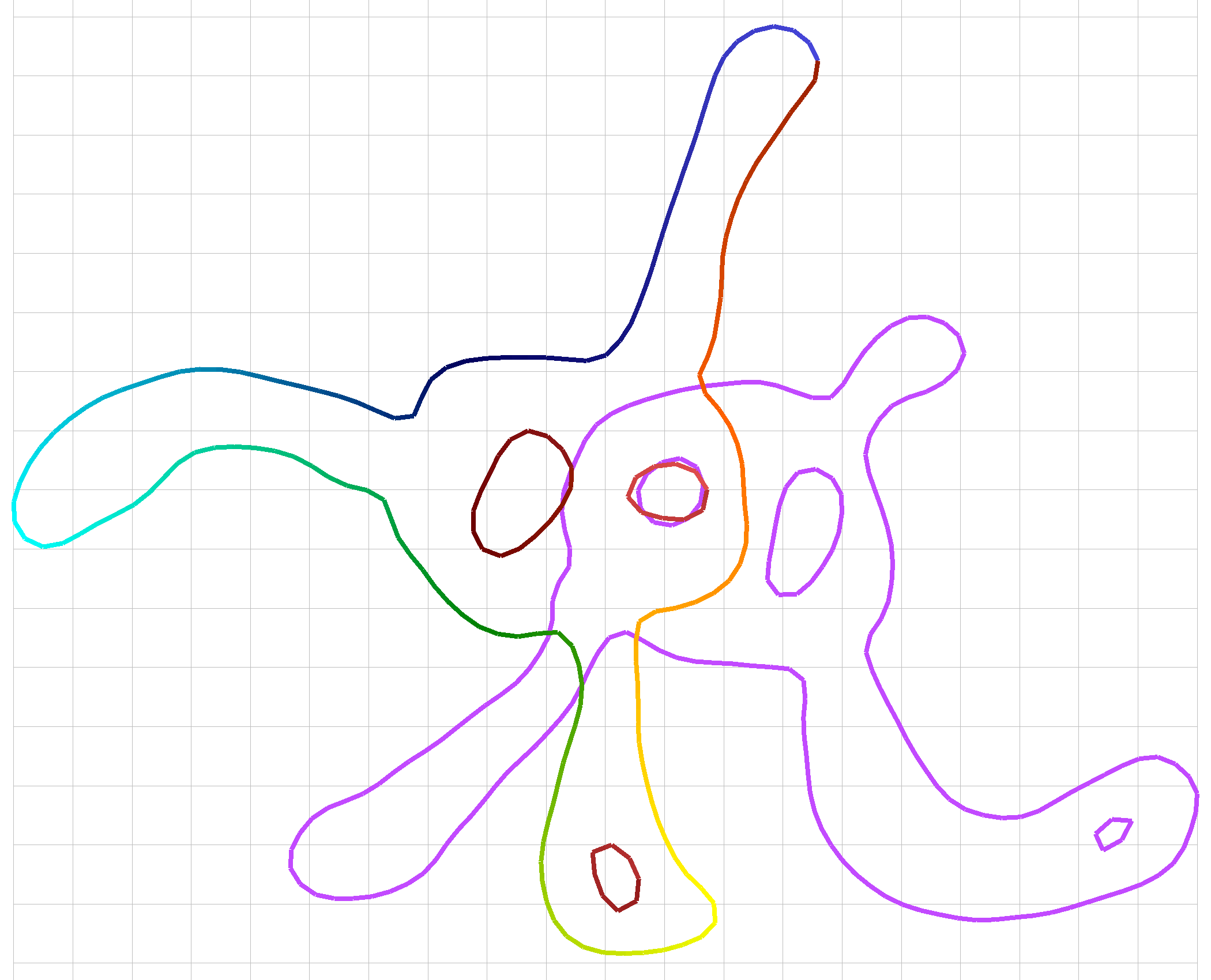}&
\includegraphics[width=.22\linewidth]{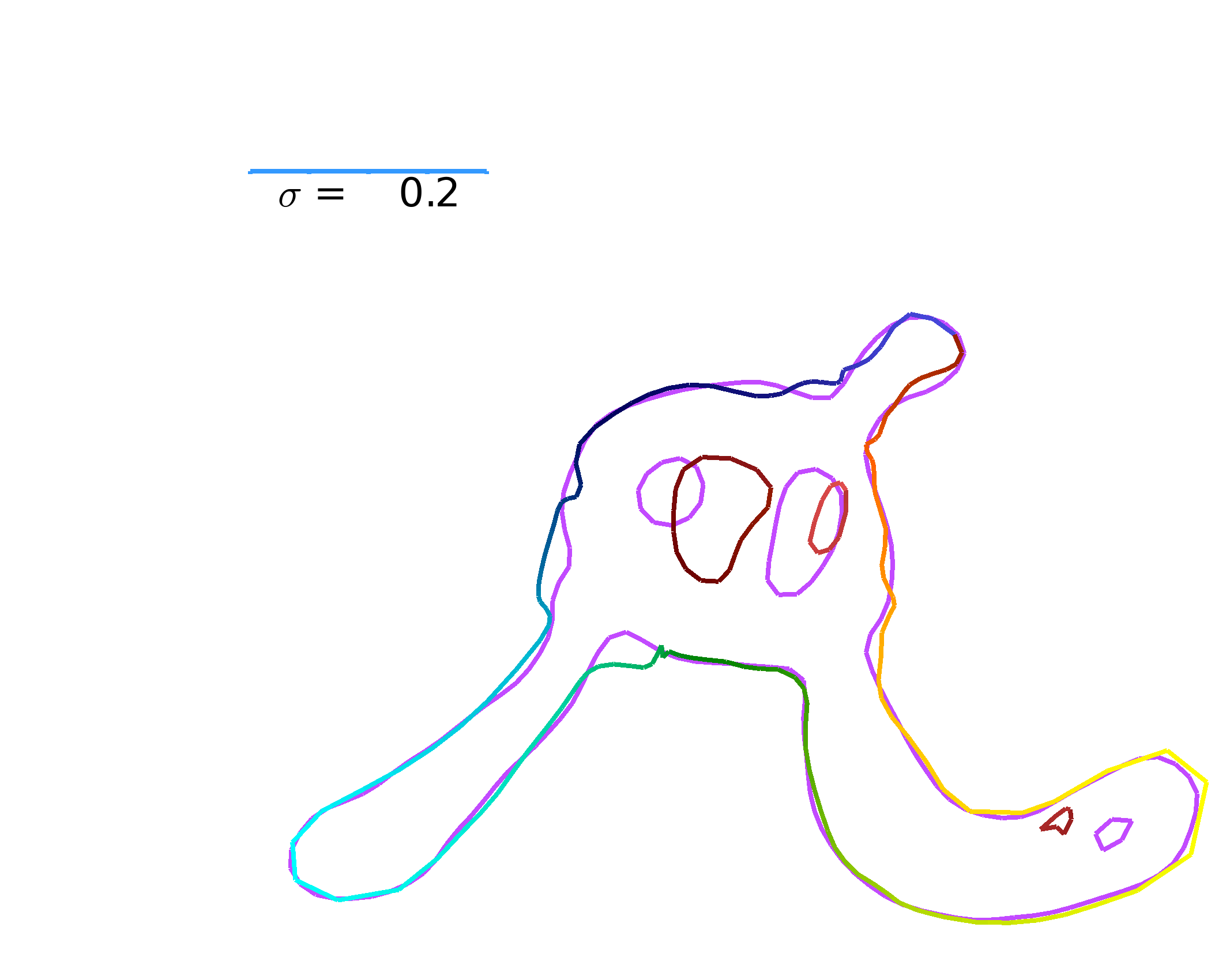}&
\includegraphics[width=.22\linewidth]{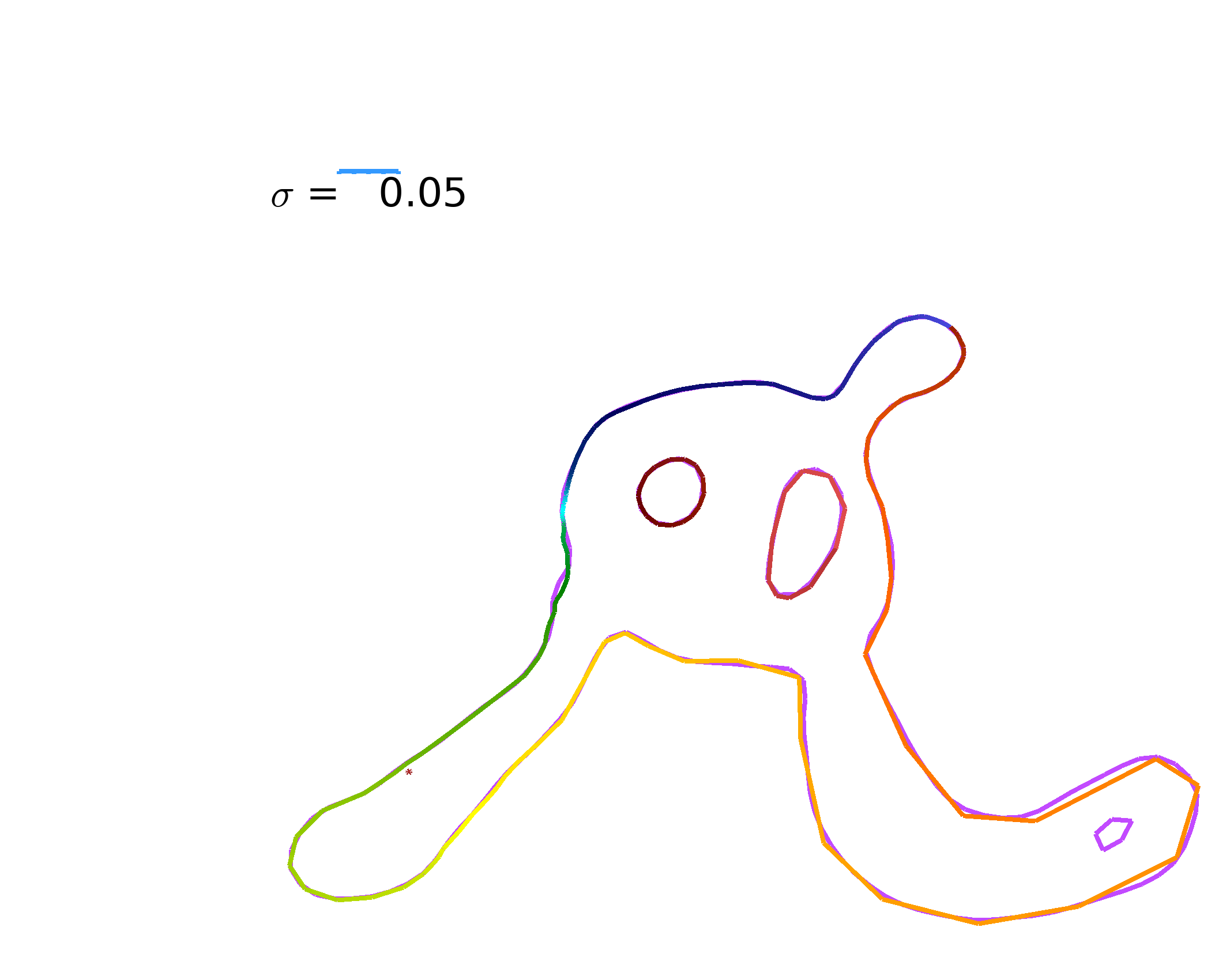}\\
(a) Dataset & (b) RKHS matching $\sigma=.2$ & (c) $\sigma=.05$
\end{tabular}
\begin{tabular}{@{}c@{}c@{}c@{}c@{}c@{}}
\includegraphics[width=.19\linewidth]{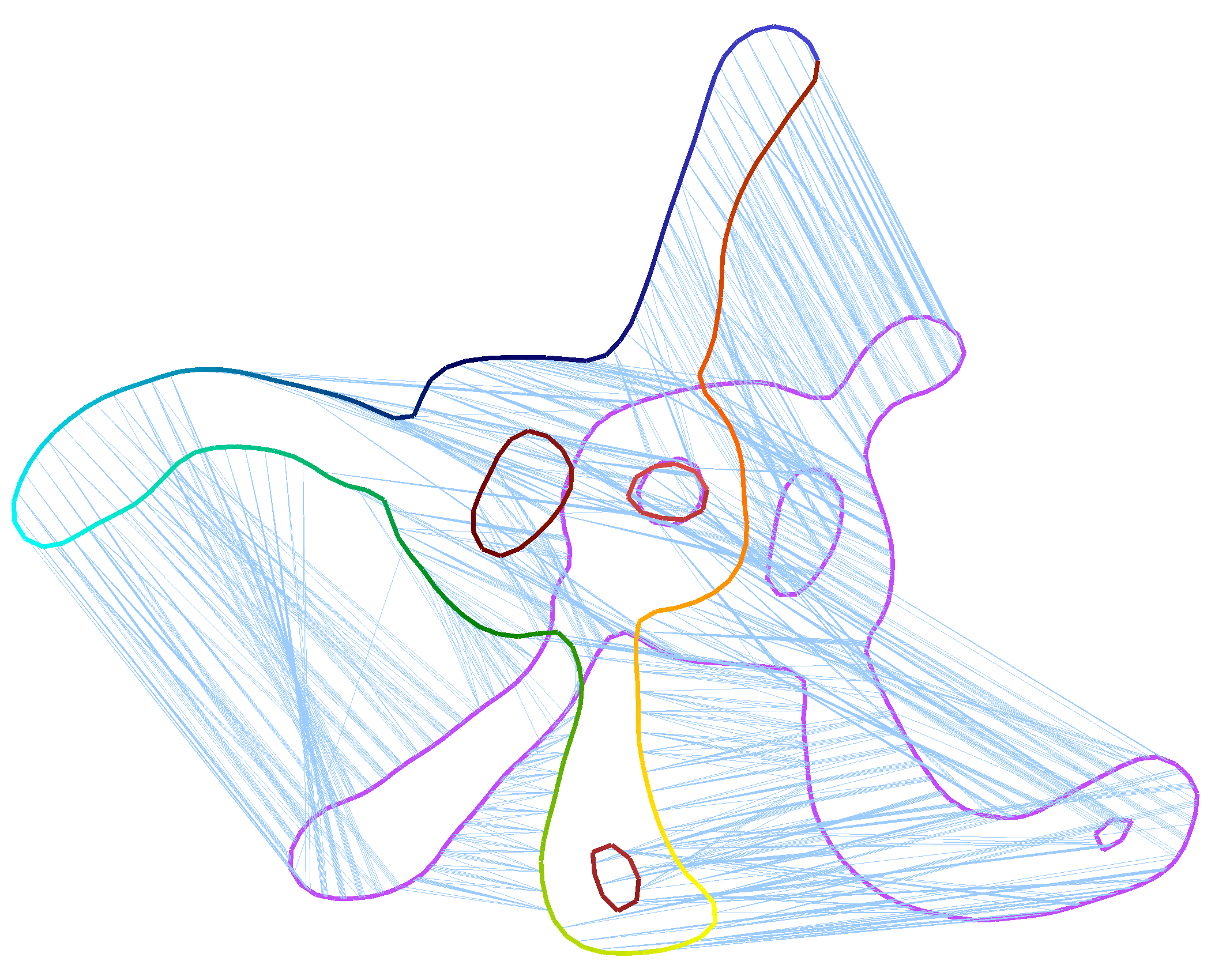}&
\includegraphics[width=.19\linewidth]{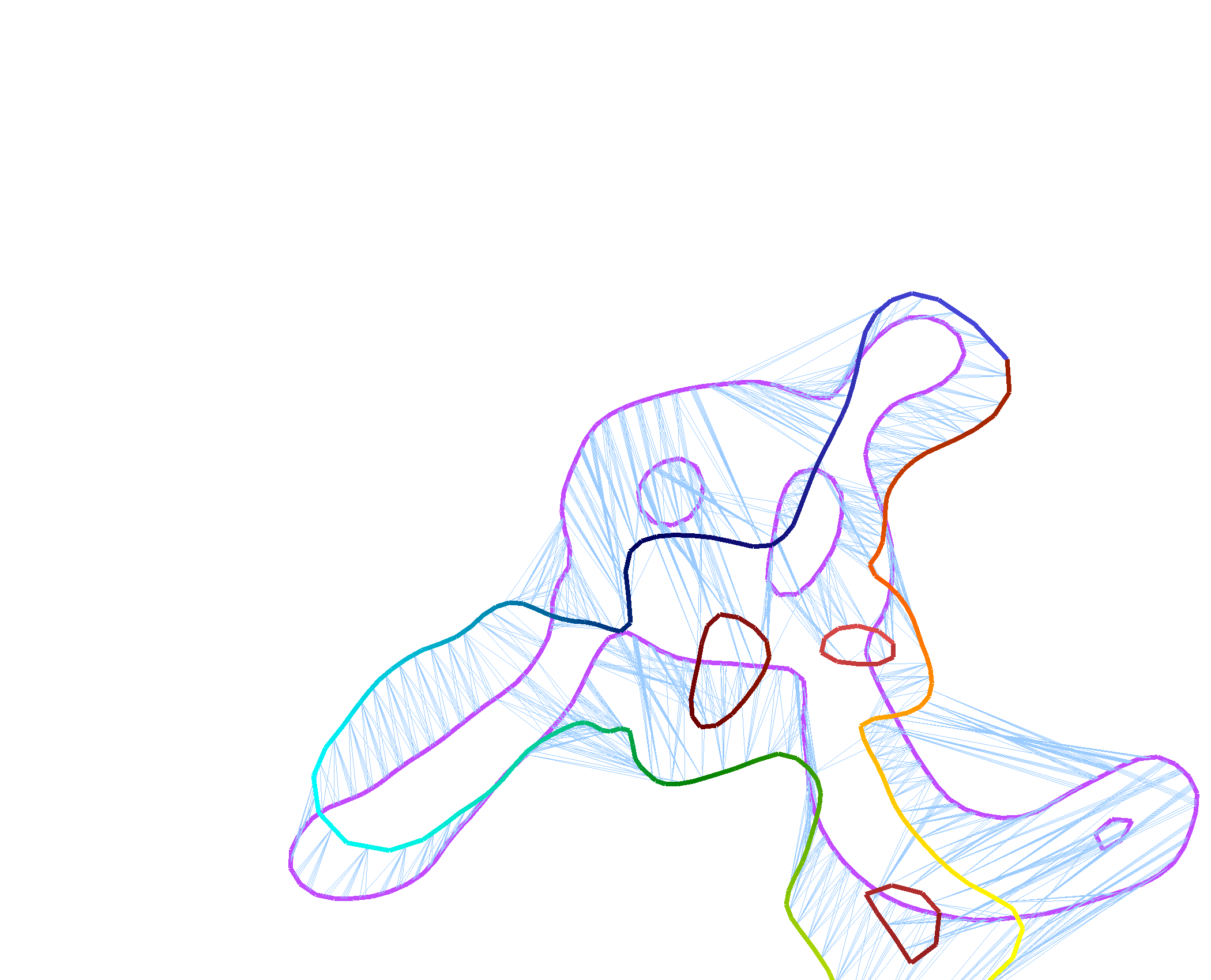}&
\includegraphics[width=.19\linewidth]{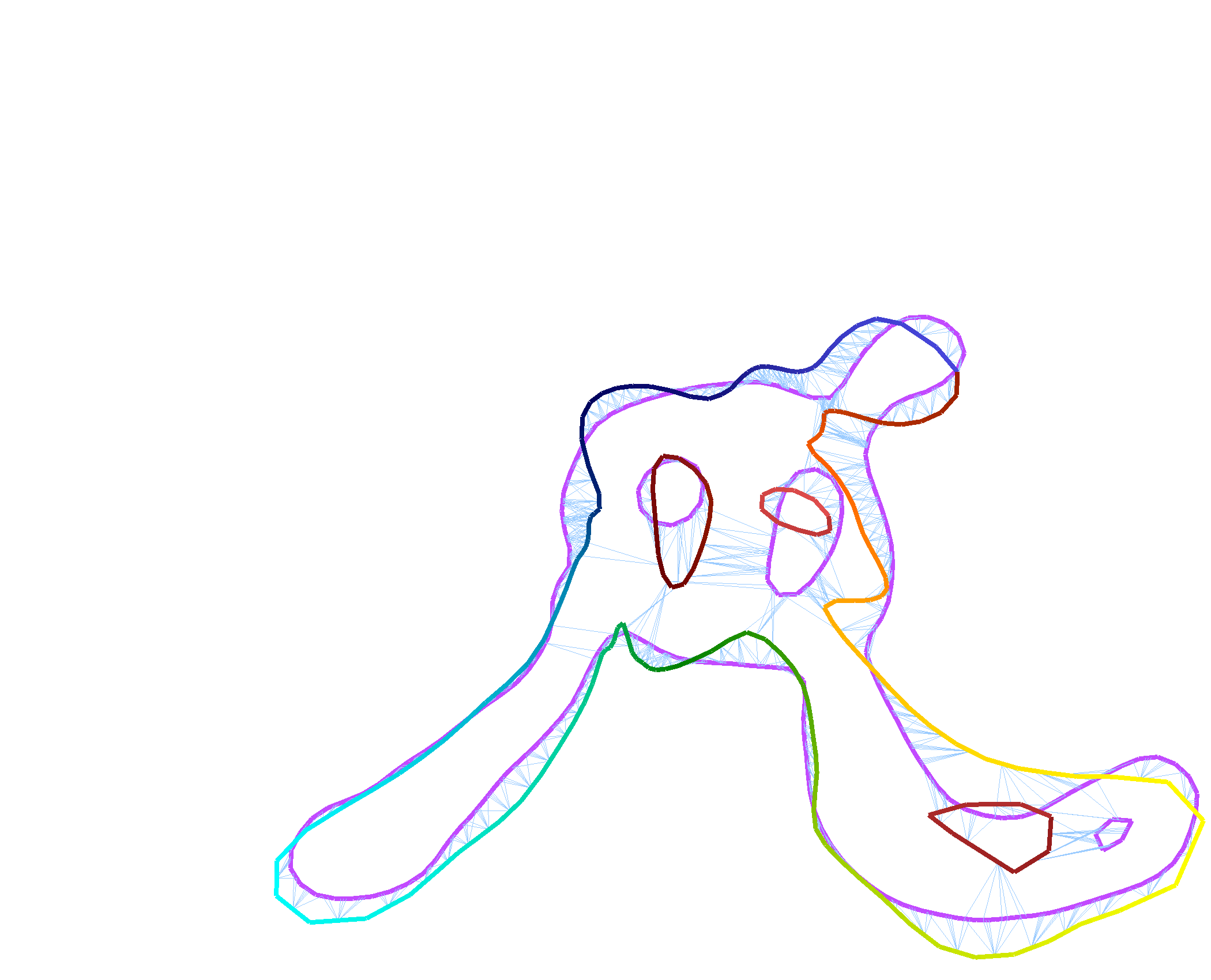}&
\includegraphics[width=.19\linewidth]{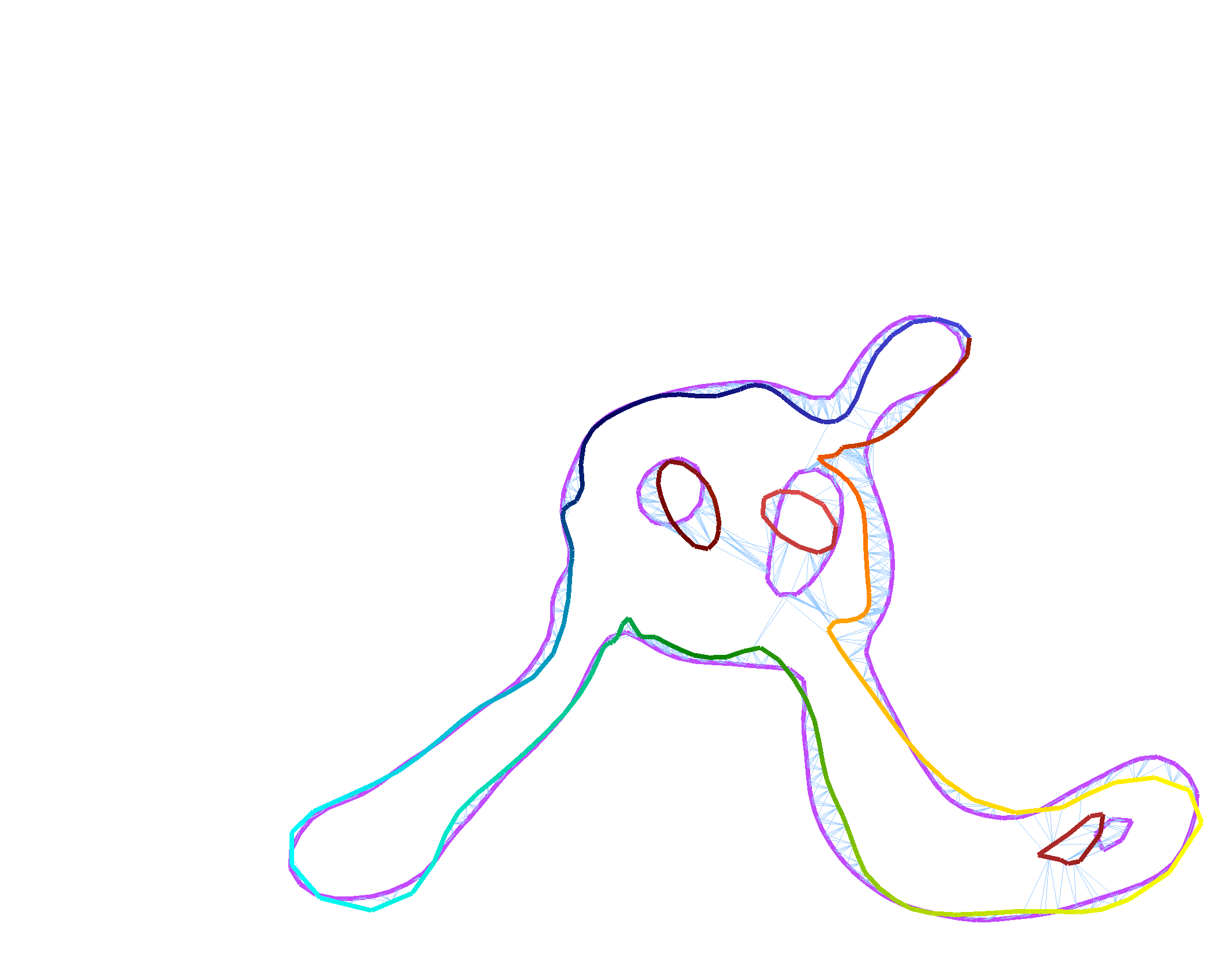}&
\includegraphics[width=.19\linewidth]{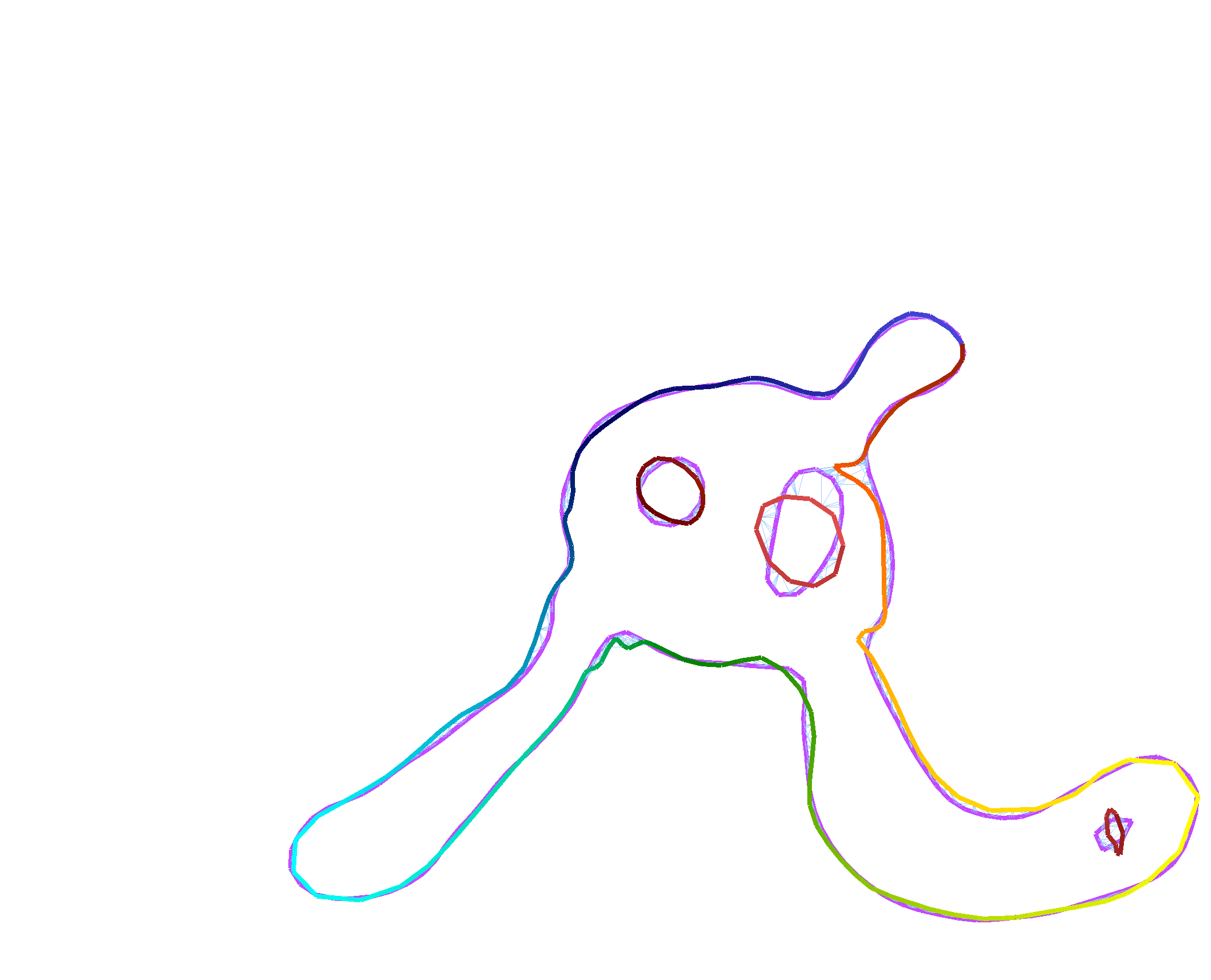} \\
$1^{\text{st}}$ fidelity computed & $5^{\text{th}}$ & 
$10^{\text{th}}$ & $20^{\text{th}}$ & $40^{\text{th}}$
\end{tabular}
\begin{tabular}{@{~~~}c@{}c@{}c@{}c@{}c@{}}
\includegraphics[width=.24\linewidth]{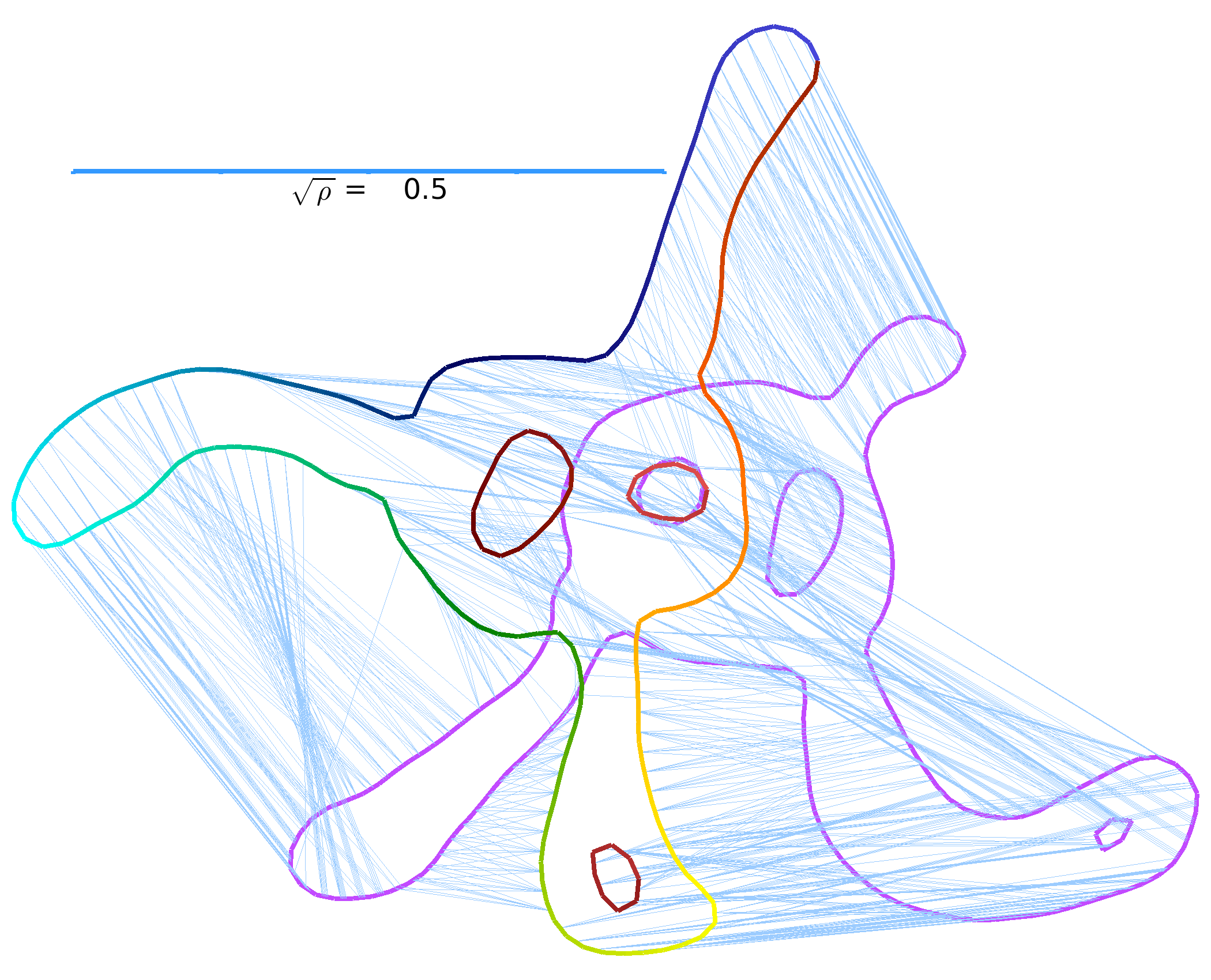} & 
\includegraphics[width=.24\linewidth]{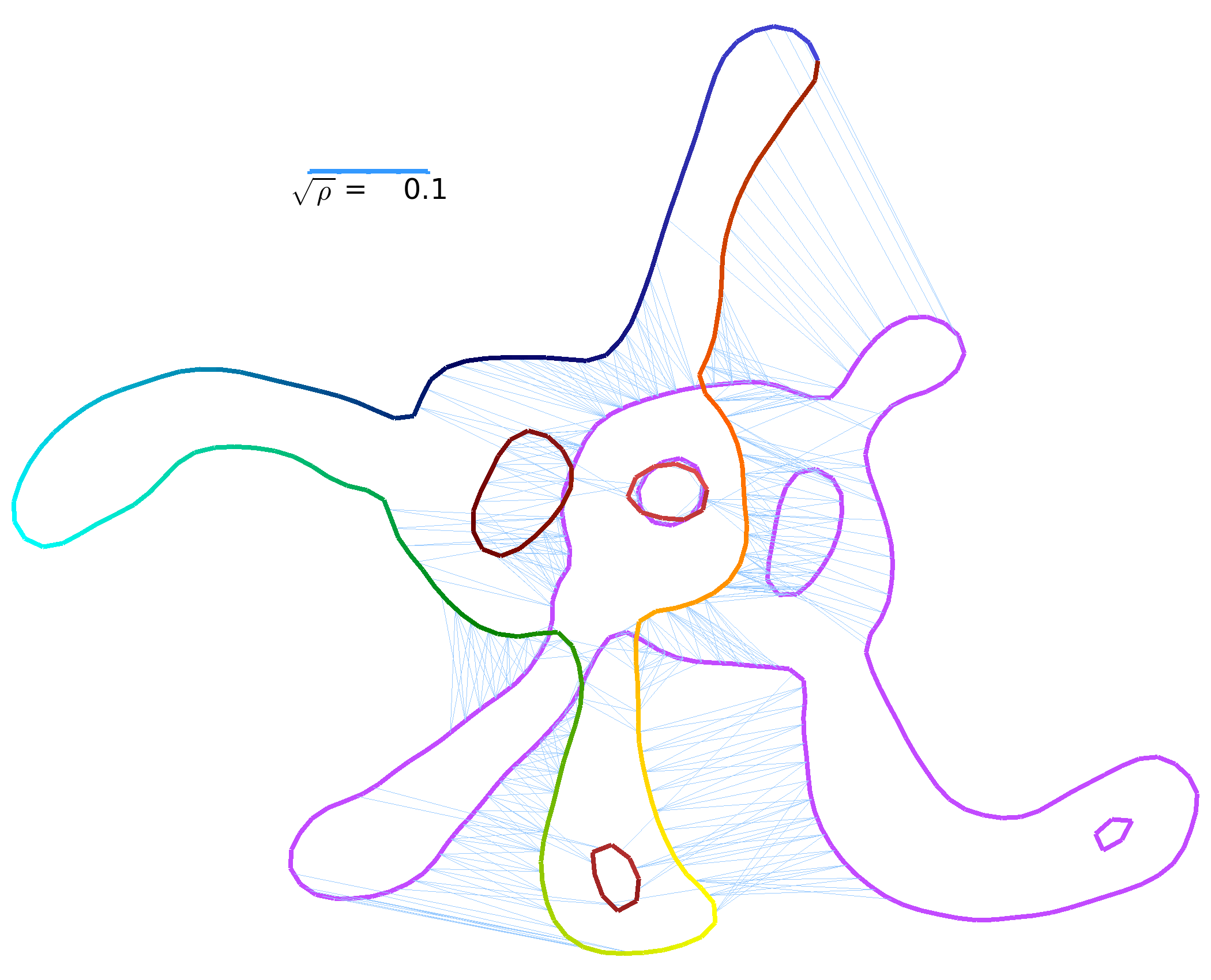} &
\includegraphics[width=.24\linewidth]{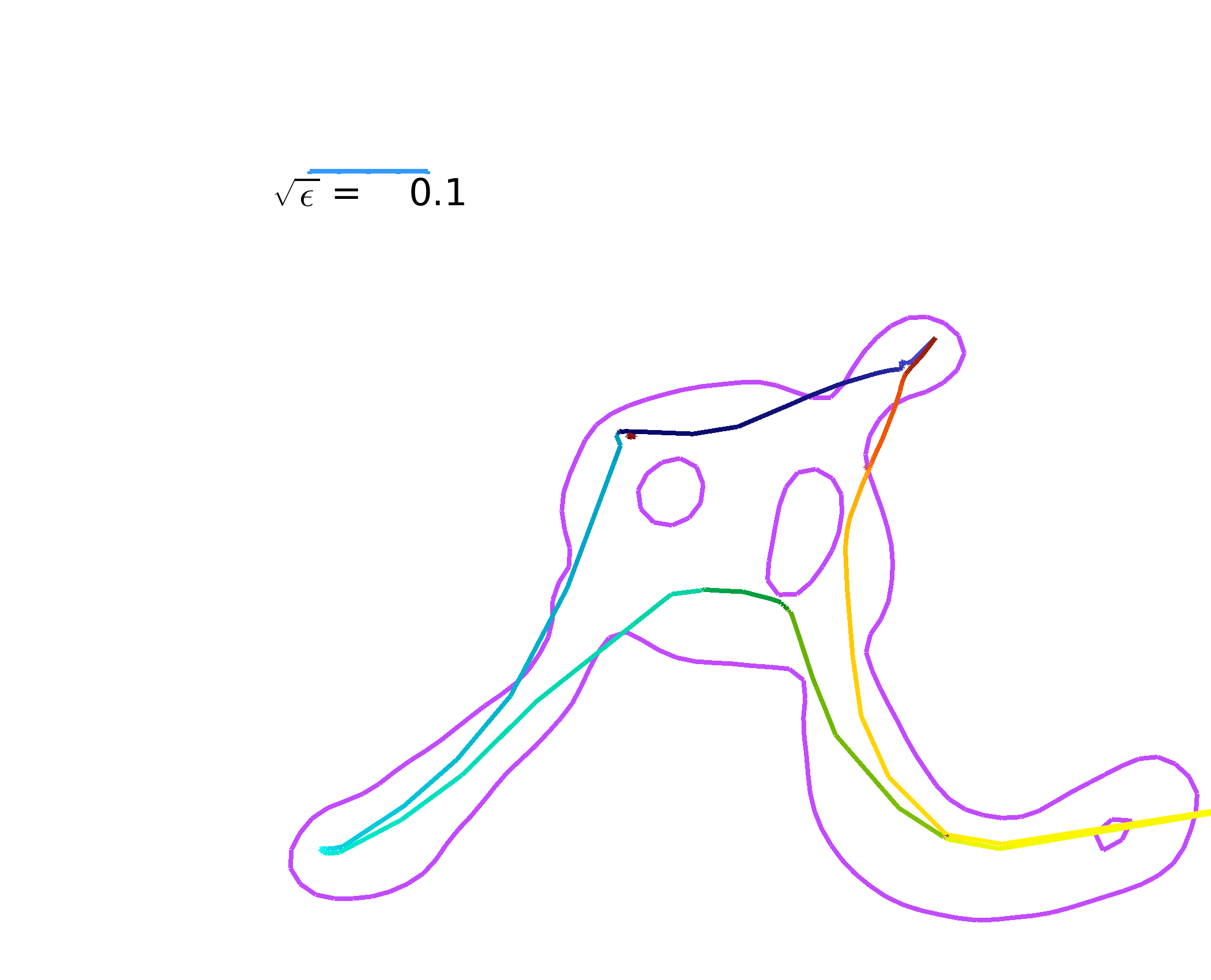} & 
\includegraphics[width=.24\linewidth]{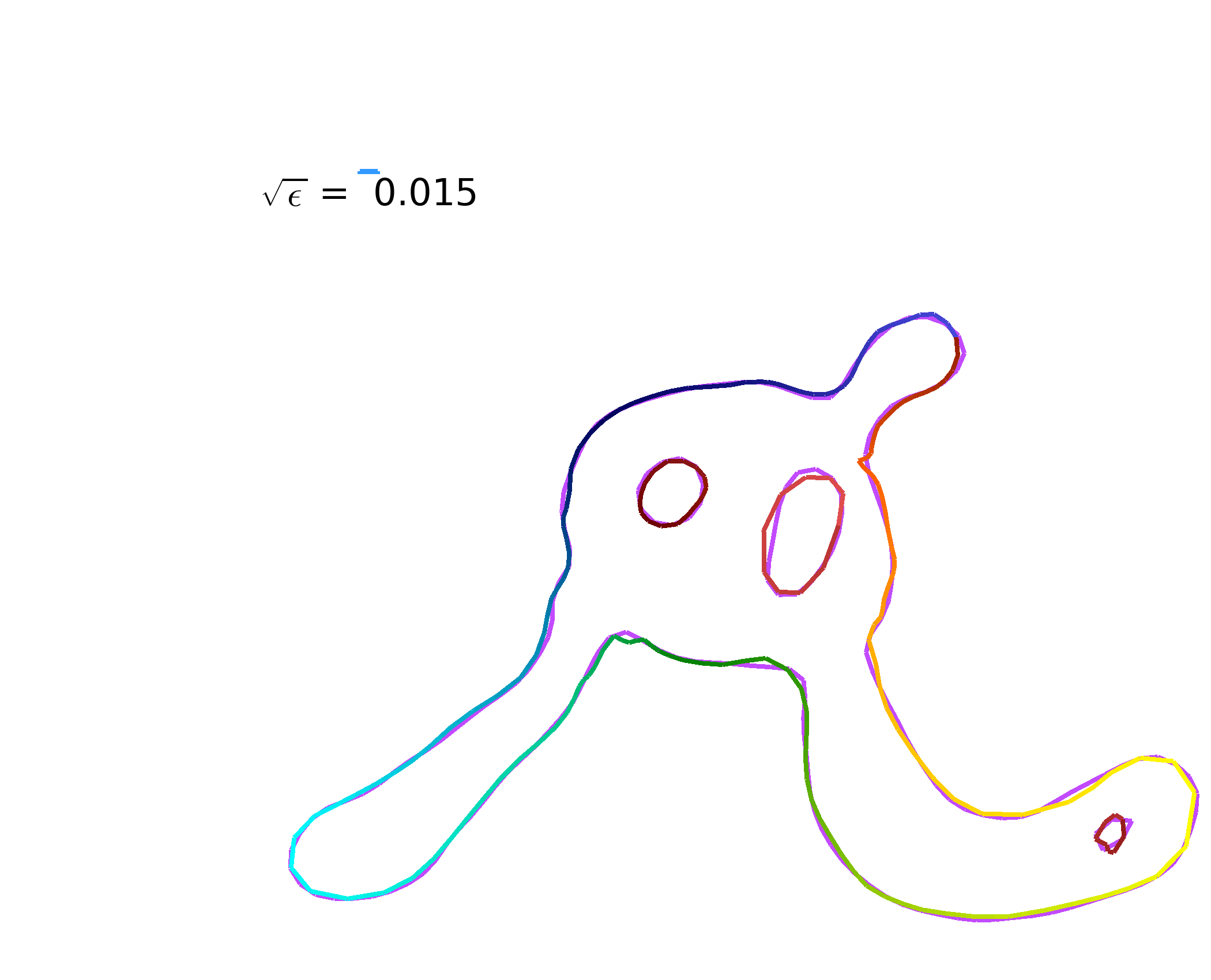}
\\
(d) $\sqrt{\rho} = .5$ & (e) $\sqrt{\rho} = .1$ & 
(f) $\sqrt{\epsilon} = .1$ & (g) $\sqrt{\epsilon} = .015$ 
\end{tabular}
	\caption{
	\textbf{First row:} presentation of a difficult registration problem.
	Even though it looks precise, (c) completely mismatches the shapes' arms
	as evidenced by the \emph{color code}.
	\textbf{Second row:} evolution of the registration algorithm (minimizing $\Ee$). 
	\textbf{Third row:} influence of $(\rho,\epsilon$). 
		\label{fig:introduction}}
\end{figure}

%
The registration is then obtained using the fidelity $\Loss=W_{\epsilon,\rho}$ in the registration problem~\eqref{eq-min-register}. 
The impact of this change (with respect to using more classical fidelity terms such as~\eqref{eq-rkhs-loss}) simply requires to input the expressions~\eqref{eq-ot-gradients} of the gradients in the chain rule~\eqref{eq-chain-rule}, which can be evaluated after running Sinkhorn algorithm~\eqref{eq-sinkhorn} to compute the optimal $\ga$ and $u$ needed in formula~(\ref{eq-ot-gradients}).
%
%
In order to get non-negative fidelities, 
one can also discard the entropy and KL divergence
from the final evaluation of the cost $\mathcal{E}$,
and compute its derivatives using an autodiff library such as Theano~\cite{theano2016}~:
this is what was used for Figure~\ref{fig:introduction}.




\mysubsection{Synthetic dataset} \label{sec:protozoa}
Figure~\ref{fig:introduction} showcases the results of our methods on a difficult 2-D curve registration problem. The first curve (rainbow colors, represented as a measure $\mu$) is deformed into the purple one (measure $\nu$).  Both curves are rescaled to fit into the unit square, the background grid of (a) is graduated every .05, and the cost function used is that of equation~\eqref{eq:cprodbis} with $\alpha=1$, $k=4$.
The diffeomorphism is computed with an LDDMM sum of Gaussian kernels, 
with $k(x,y) = 1.\,\exp(-\norm{x-y}^2/(2\cdot.025^2)) + .75 \,\exp(-\norm{x-y}^2/(2 \cdot .15^2))$.

\myparagraph{RKHS fidelity: first row (b),(c)}
%
(b) and (c) have been computed using a \emph{kernel-varifold} fidelity, 
with a spatial Gaussian kernel of deviation $\sigma$ and an acute angular selectivity 
in $\cos^4(\theta)$ -- with $\theta$ the angle between two normal directions.
As shown in (b), such an RKHS fidelity method performs well with a large bandwidth $\sigma$. Unfortunately, trying to increase the precision by lowering the value of $\sigma$ leads to the creation of 
undesirable local minima. In the eventual registration (c) the arms are not transported, but shrinked/expanded, as indicated by the \emph{color code}.  Classical workarounds involve decreasing $\sigma$ during the shape optimization in a coarse-to-fine scheme, which requires a delicate parameters tuning. The main contribution of this paper is that it provides a principled solution to this engineering problem, which is \emph{independent} of the underlying optimization/gradient descent toolbox, and can be adapted to any non-local fidelity term.

\myparagraph{OT fidelity: second row}
In sharp contrast with this observed behavior of RKHS fidelity terms, the OT data attachment term overcomes local minima through the computation of a \emph{global} transport plan, displayed in light blue. 
Note that since the cost function used in this section is \emph{quadratic}, $\rho$ and $\epsilon$ should be interpreted as squared distances, and we used here $\sqrt{\epsilon} = .015$, $\sqrt{\rho}=.5$. The final matching is displayed (g).

\myparagraph{Influence of $\rho$, third row (d),(e)}
Here, we used $\sqrt{\epsilon} = .03$.
The value of  $\sqrt{\rho}$ acts as a ``cutoff scale'', above which OT fidelity tends to favour mass destruction/creation over transport. 
This result in a ``partial'' and localized transport plan, which is useful when dealing with outliers, large mass discrepancies which should not be explained through transport.

%

\myparagraph{Influence of $\epsilon$, third row (f),(g)}
Here, $\sqrt{\rho} = .5$.
$\sqrt{\epsilon}$ should be understood as a diffusion, blurring scale on the optimal transport.
The resulting matching can therefore only capture structure up to a scale $\sqrt{\epsilon}$~:
in (f), the ``skeleton'' mean axis of the shape.

\myparagraph{Computational cost}
The number of steps needed to compute a transport plan roughly scales like $O(\rho/\epsilon)$.
In this experiment, an evaluation of the fidelity term and its gradient was 100 (resp. 1000) times as long to compute as a RKHS loss of the form~\eqref{eq-rkhs-loss}, for $\sqrt{\epsilon} = .1$ (resp $.015$). It thus has roughly the same complexity (resp. one order of magnitude slower) than evaluating the LDDMM diffeomorphism itself.
As shown in the second row of Figure~\ref{fig:introduction}, an efficient
optimization routine may only require to evaluate the OT fidelity a handful of times to be driven to an appropriate rough deformation.
Although not used here, a heuristic to drastically reduce the computational workload is a two-step scheme: first, use OT with a large $\epsilon$ to find the good basin of attraction ; then, use a fast non-local fidelity (e.g.~\eqref{eq-rkhs-loss} with small $\sigma$) to increase precision.

\begin{figure}[!b]
\centering
\begin{tabular}{@{}c@{}c@{}c@{}}
\includegraphics[width=.24\linewidth]{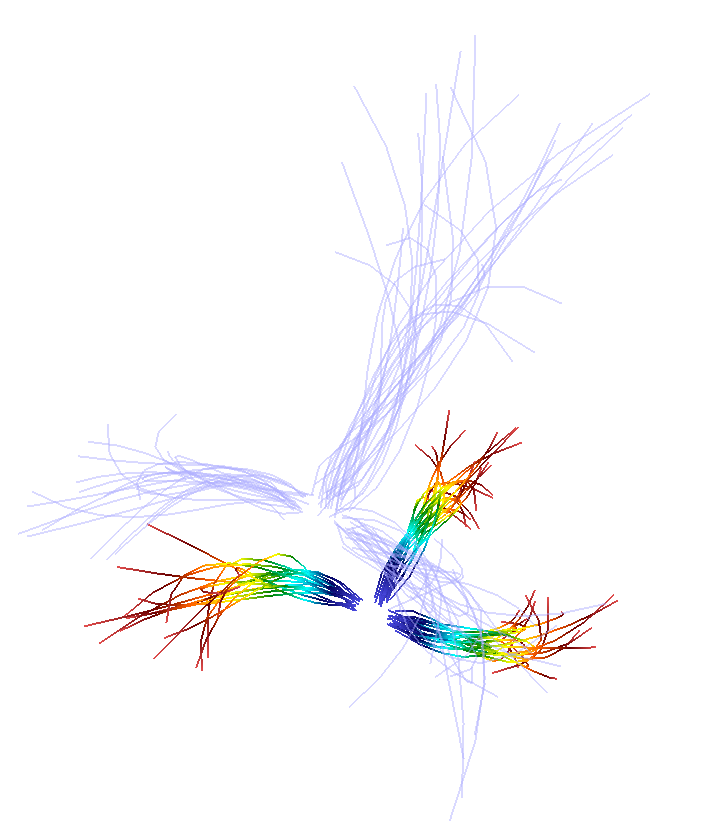}&
\includegraphics[width=.24\linewidth]{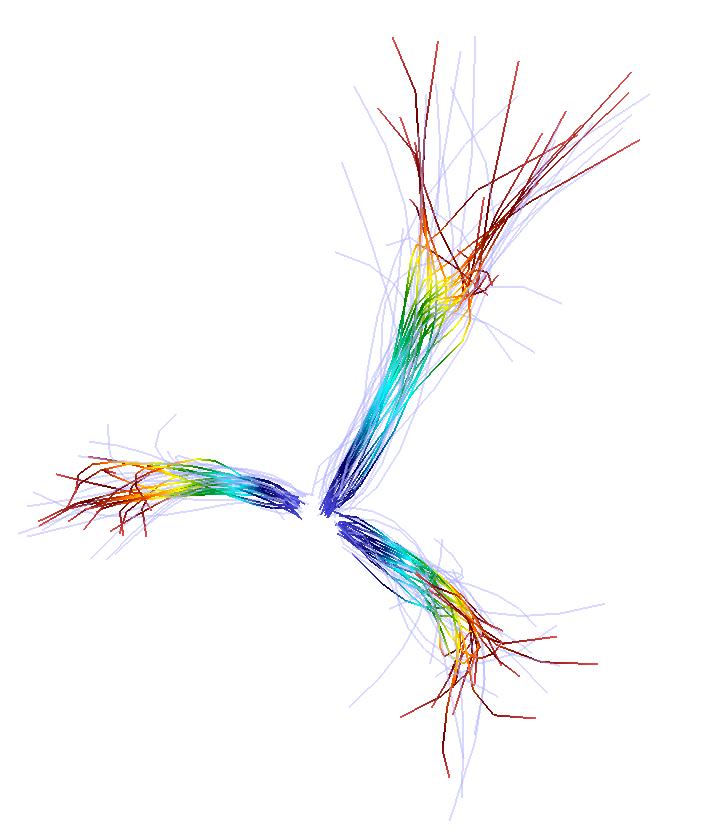}&
\includegraphics[width=.24\linewidth]{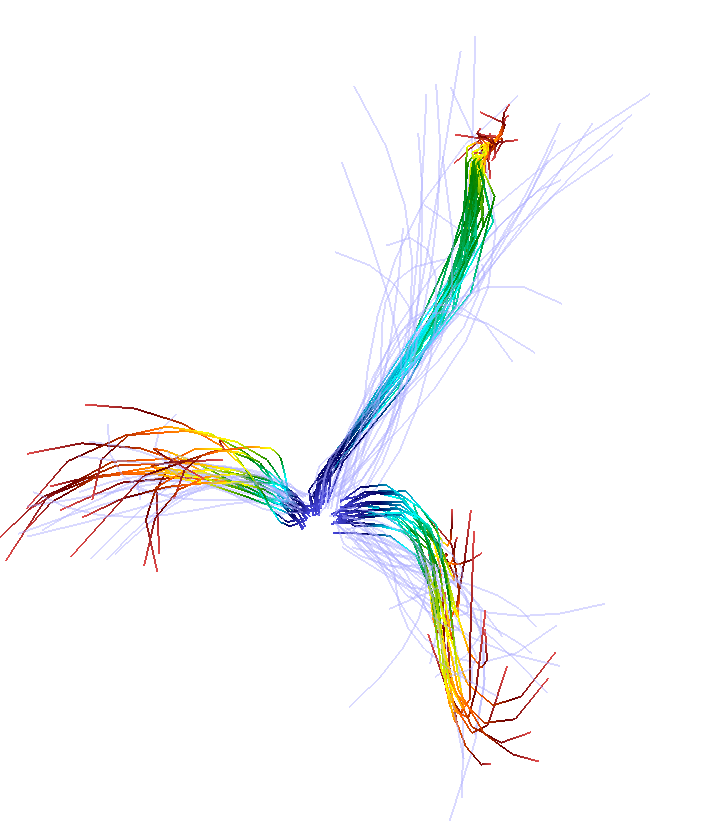} \\[-2mm]
(a) Dataset & (b) OT fidelity & (c) RKHS fidelity
\end{tabular}\vspace{-1mm}
\begin{tabular}{@{}c@{}c@{}c@{}c@{}}
\includegraphics[width=.2\linewidth]{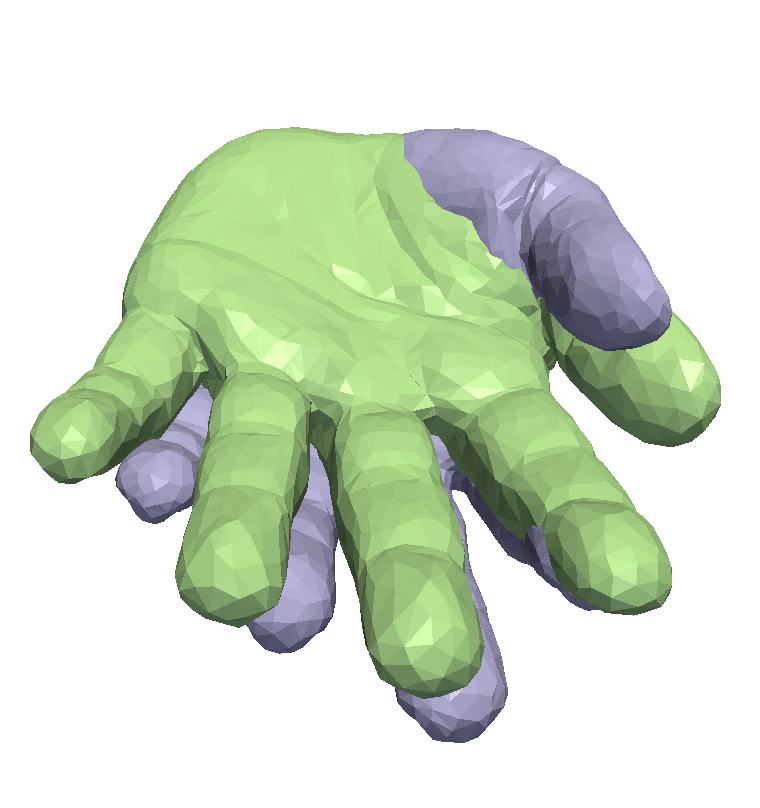}&
\includegraphics[width=.2\linewidth]{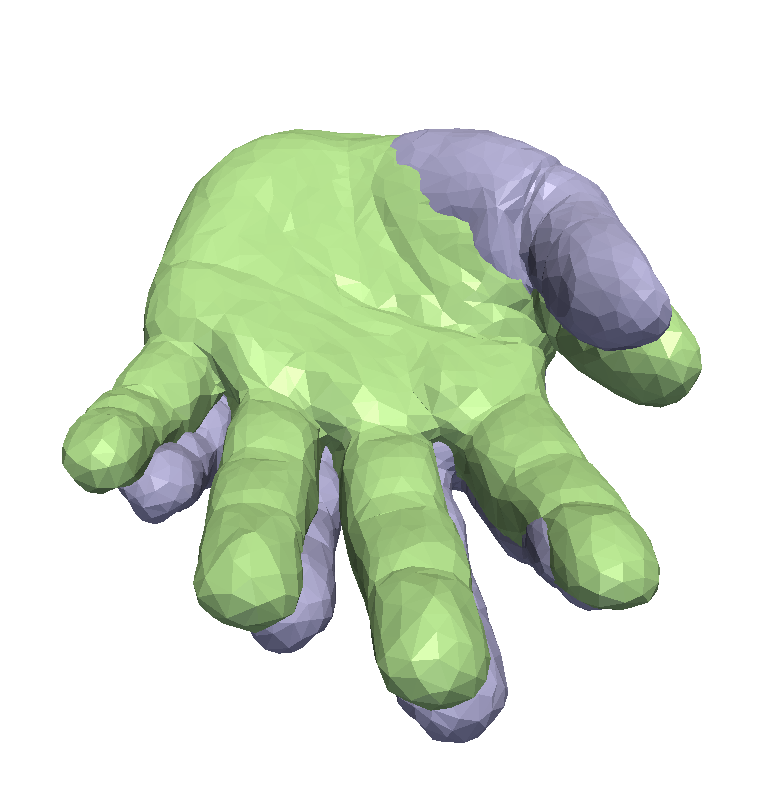}&
\includegraphics[width=.2\linewidth]{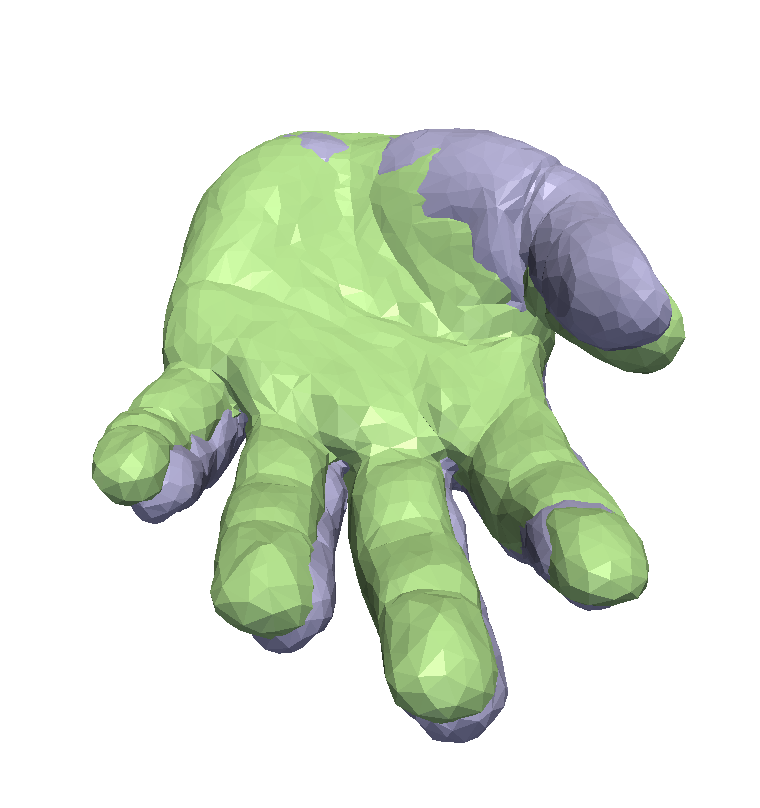}&
\includegraphics[width=.2\linewidth]{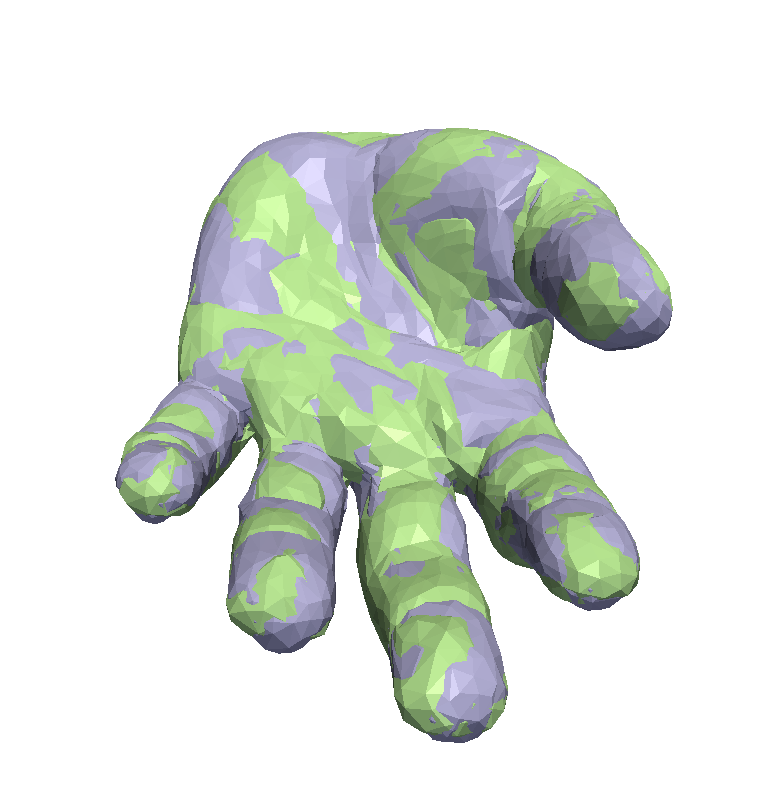}\\[-2mm]
$t=0$ & $t=1/3$ & $t=2/3$ & $t=1$
\end{tabular}\vspace{-3mm}
\caption{%
\textbf{First row:} Matching of fibres bundles.
\textbf{Second row:} Matching of two hand surfaces using a balanced OT fidelity. Target is in purple.
}
\label{fig:fb}\label{fig:hands}
\end{figure}

\mysubsection{Fibres bundles dataset}
The numerical experiment presented Figure \ref{fig:fb} illustrates the problem of registration of fibres bundles in 3d. It is often difficult to compute convincing registration of fibers bundles as the ends of the fibres are in practice difficult to align. This toy example may be seen as a very simple prototype of white matter bundles common in brain imaging. Currents-based distance together with a LDDMM framework were already used to analyze this kind of data, see e.g; \cite{gori14}. 

The source and target shape have 3 bundles containing around 20 fibers each. The diameter of the dataset is normalized to fit in a box of size 1. The cost function used for the OT fidelity is \eqref{eq:cprod} with the orientation dependant distance between normals. We use the unbalanced framework with $\sqrt{\rho}=1$ and $\sqrt{\varepsilon}=0.07$. Using this OT fidelity with LDDMM allows to recover the shape of the target bundles (see Figure \ref{fig:fb} first row) whereas the RKHS fidelity based registration (we use a Gaussian kernel width $\sigma = 0.8$) converges toward a poor local minimum.

\mysubsection{Hands dataset}
OT fidelity may be used with large datasets thanks to an efficient implementation of Sinkhorn iterations. The two hand shape surfaces of Figure~\ref{fig:hands} contain more than 5000 triangles. The registration takes less than 1 hour on a GPU. This numerical experiment shows that OT fidelity may be used to register surfaces with features at different scales. 

\section*{Conclusion}

In this article, we have shown that optimal transport fidelity leads to more robust and simpler diffeomorphic registration, avoiding poor local minima.
%
Thanks to the fast Sinkhorn algorithm, this versatile tool has proven to be usable and scalable on real data and we illustrated its efficiency on curves, surfaces and fibres bundles. We plan to extend it to segmented volumetric image data.

{\tiny{
\bibliographystyle{splncs03}      
\bibliography{biblio}    
}}
\end{document}